# Researches on the Twin Prime Problem

**Zhang Baoshan**[①]

(School of Science, Nanjing Audit University, Nanjing, China, Zip code: 211815)

Email: bszhang@nau.edu.cn bszhang8@163.com

**Abstract**: Twin prime number problem is mainly the structure of the twin prime numbers and whether there are infinitely many prime twins group. In this paper, by constructing a special cluster number set, proves that the number of set number of the first n columns set the intersection of the minimum number of $q$ is decision of the prime twins ($q$, $q$ + 2), and the minimum number of series is divergent. Prime twins so thoroughly proved there must be infinitely many groups of twin prime conjecture.

## §1 Definition and Lemmas

It is well known that natural number set $1,2,3,\cdots,n,\cdots$ can be classified according to different characteristics, thus get the new data set, such as the even, odd number set, a prime number sets, etc. For the convenience of first presents the definition of prime numbers and **composite number**s.

**Definition**1. If natural number $n$ to remove 1 and itself without the natural factor, $n$ is called a **prime number**; Other than 1 and a prime natural number $n$ called a **composite number**.

**Definition**2. If $q$, $p$ are a pair of prime numbers, and $|p-q| \le 2$, It says they are twin prime number.

Obviously, $2,3$ are the smallest prime twins and other prime twins are like $q$, $q+2$ the odd prime numbers.

Usually, using $N$ represented by the natural number set, $P$ or $\{p_k\}$ said collection of prime number, and $p_k$ is according to the natural order of the first a

---

[①]Baoshan Zhang. male, 1959 born, PhD, Professor, Native of Fengxian County in jiangsu province of China. Mainly engaged in the study of differential equations, computer symbolic operation and so on, also has research interest in number theory, published over 70 papers, published by 2 books. Email: bszhang@nau.edu.cn bszhang8@163.com

prime number. Here are some common results of elementary number theory [1 ~ 2].

**Lemma 1** There are an infinite number of prime set, $\{ p_k \}$is prime sequence strictly monotone increasing sequence of $N$ .

**Lemma 2** Let $p$ be a prime, for any natural number $m$ , Then $m, p$ is relatively prime or $m$ is divisible by $p$ , that is

$$(m, p) = 1 \text{ or } p \mid m. \qquad (1.1)$$

**Lemma 3**（**Twin Prime Conjecture** [1~3]）There is infinitely twin prime number groups such as $q$ , $q+2$ .

Lemma 1 ~ 2 is well known, and Lemma 3 (twin prime conjecture) is now a famous unsolved problems in number theory. Mathematician Hilbert transform the guess to famous international congress of mathematicians in the report in 1900, it ranked eighth in the 23 "Hilbert problem", can be described as "there are infinitely many prime twins" [4 ~ 20].Recent work in Zhang Yitang and others greatly promote the study of the problem [21 ~ 24], but there are still some distance away from the final solution.

In this paper, by constructing a special cluster number set, proves that the number of set number of the first n columns set the intersection of the minimum number of $q$ is decision of the prime twins ($q$, $q$ + 2), and the minimum number of series is divergent. Prime twins so thoroughly proved there must be infinitely many groups of twin prime conjecture.

## § 2 The structure of the derivative of twin prime number set sequence

To explore the formation of twin prime numbers and limitless problems, explore the effective ways to prove Lemma 3, we have to construct the following set sequence $\{S_n\}$ .

$$S_n = \{ a \in N \mid a = 2k + 1 > n, (a, n) = (a + 2, n) = 1, n \in P \} , \qquad (2.1)$$

For example,

$S_1 = \{ 2k+1 \mid \forall k \in N \}$,

$S_2 = \{ a > 2 \mid (a,2) = (a+2,2) = 1 \} = \{ 3,5,7,9,\cdots,2n+1,\cdots \}$,

$S_3 = \{ a > 3 \mid (a,3) = (a+2,3) = 1 \} = \{ 5,11,17,\cdots,\cdots \}$,

$\cdots,\cdots,\cdots,$

$$S_{p_n} = \{ a = 2k+1 > p_n \mid (a, p_n) = (a+2, p_n) = 1 \}, \cdots\cdots, \qquad (2.2)$$

And here $p_n$ said the *n-th* prime. Described above sequence's structure is $S_1$ form of the all odd numbers ; $S_2$ is the all odd numbers more than 2, and $\min S_2 = 3$; $S_3$ is the set of odd number $a$, which is more than 3 and $a,\ a+2$ are relatively prime with $p_2 = 3$, and also $\min S_3 = 5$; $S_5$ is the set of odd number $a$, which is more than 3 and $a,\ a+2$ are relatively prime with $p_3 = 5$, and $\min S_5 = 7$; And so on , $S_{p_n}$ is the set of odd number $a$, which is more than 3 and $a,\ a+2$ are relatively prime with $p_n$, and then $\min S_{p_n} = p_n + 2$.

By (2.1), (2.2) , we get the set sequence

$$\{ S_{p_n} \}_{n=1}^{\infty} = \{ S_{p_1}, S_{p_2}, \cdots, S_{p_n}, \cdots \}. \qquad (2.3)$$

It has the following features.

**Theorem 1**. The set sequence $\{ S_{p_n} \}_{n=1}^{\infty}$ has property: 1) $\Phi \neq S_{p_n} \subset S_1$ is true for all $n$ ; 2) For any natural number M, the $\bigcap_{n=1}^{M} S_{p_n} \neq \Phi$ will be true.

**Proof.** By formula (2.1) and the structure of the process about $S_2, S_3, S_5, \cdots, S_{p_n}, \cdots$,we have $\Phi \neq S_{p_n} \subset S_1$, That is the conclusion 1) is true. In fact, it is note that $p_n$ is a prime, there would be $(q,\ p_n) = 1$, $(q+2,\ p_n) = 1$ for any odd prime number $q$ in the interval $[p_n, 2p_n]$, and then $q \in S_{p_n}$.

Let's certificate the conclusion 2), for any natural number M, choosing

$$q = \prod_{n=2}^{M} p_n = p_2 p_3 \cdots p_M \quad , \qquad (2.4)$$

and by using of prime properties, it easy to know

$$(q+2,\ p_j)=1,\quad (q+4,\ p_j)=1,\quad j=2,3,\cdots M.$$

Therefore, we have

$$q+2\in S_{p_j},\quad j=2,3,\cdots M \tag{2.5}$$

That is $q+2\in\bigcap_{n=2}^{M}S_{p_n}$, and so $\bigcap_{n=2}^{M}S_{p_n}\neq\Phi$. This is the proof of Theorem 1.

According to the conclusion of Theorem 1, we can construct a sequence $\{q_M\}_{M=2}^{\infty}$ by through the collection of the minimum value to $\bigcap_{n=2}^{M}S_{p_n}$ as below.

**Theorem 2.** The set sequence $\{S_{p_n}\}_{n=1}^{\infty}$ has property: For any natural number M, Considering the minimum value of $\bigcap_{n=2}^{M}S_{p_n}$ :

$$q_M=\min\left(\bigcap_{n=2}^{M}S_{p_n}\right). \tag{2.6}$$

If $q_M\in[p_M,2p_M]$, then $q_M,\ q_M+2$ is a pair of the prime twins.

**Proof.** From (2.6), we have, $q_M\in S_{p_n}$, $(n=2,3,\cdots,M)$, thus

$$(q_M,\ p_n)=1,\quad (q_M+2,\ p_n)=1,\quad (n=2,3,\cdots,M). \tag{2.7}$$

By Lemma 2, $q_M,\ q_M+2$ will be a pair of the prime twins in case of $q_M\in[p_M,2p_M]$. And so, $q_M,\ q_M+2$ are prime twins. This is the proof of Theorem 2.

The following Mathematica program can achieve certain values under the calculation problem (2.6) and validation for some natural number M :

```
For[n = 1, n <= 1000, n++, q[n] = {};
  For[m = Prime[n] + 1, m <= Prime[n] + 3000, m++,
  If[OddQ[m] && GCD[m, Prime[n]] == 1 && GCD[m + 2, Prime[n]] == 1,
      q[n] = Join[q[n], {m}]]]]
gp[1] = q[1];
For[n = 2, n <= 400, n++, Print[gp[n] = Intersection[gp[n - 1], q[n]]]]
h = {};
For[n = 2, n <= 400, n++,
```

```
If[PrimeQ[gp[n][[1]]] && PrimeQ[gp[n][[1]] + 2], h = Join[h, {n}]]]
h
Length[h]
t = Table[{Min[gp[n]], Min[gp[n]] + 2}, {n, 2,400}];
s = Union[t]
PrimeQ[s]
Length[%]
```

Running above Mathematica program (Note to third from bottom line!)Will produce 78 of the prime twins as below：

{5, 7}, {11, 13}, {17, 19}, {29, 31}, {41, 43}, {59, 61}, {71, 73}, {101, 103}, {107, 109}, {137, 139}, {149, 151}, {179, 181}, {191, 193}, {197, 199}, {227, 229}, {239, 241}, {269, 271}, {281, 283}, {311, 313}, {347, 349}, {419, 421}, {431, 433}, {461, 463}, {521, 523}, {569, 571}, {599, 601}, {617, 619}, {641, 643}, {659, 661}, {809, 811}, {821, 823}, {827, 829}, {857, 859}, {881, 883}, {1019, 1021}, {1031, 1033}, {1049, 1051}, {1061, 1063}, {1091, 1093}, {1151, 1153}, {1229, 1231}, {1277, 1279}, {1289, 1291}, {1301, 1303}, {1319, 1321}, {1427, 1429}, {1451, 1453}, {1481, 1483}, {1487, 1489}, {1607, 1609}, {1619, 1621}, {1667, 1669}, {1697, 1699}, {1721, 1723}, {1787, 1789}, {1871, 1873}, {1877, 1879}, {1931, 1933}, {1949, 1951}, {1997, 1999}, {2027, 2029}, {2081, 2083}, {2087, 2089}, {2111, 2113}, {2129, 2131}, {2141,2143}, {2237, 2239}, {2267, 2269}, {2309, 2311}, {2339, 2341}, {2381, 2383}, {2549, 2551}, {2591, 2593}, {2657, 2659}, {2687, 2689}, {2711, 2713}, {2729, 2731}, {2789, 2791}.

It is worth pointing out, theorem 2 shows that the fact that it is set in proportion to the number (2.1) of tectonic sequence $\{q_n\}_{n=2}^{\infty}$, based on the $\{S_{p_n}\}_{n=1}^{\infty}$ are derived from the type (2.6). Corresponding sequence is derived a two-dimensional array:

$$\{(q_n,\ q_n+2)\}_{n=2}^{\infty},\quad q_n = \min\left(\bigcap_{j=2}^{n} S_{p_j}\right), \qquad (2.8)$$

Here is $q_n > p_n$. By theorem 2, if it satisfies $q_n < 2p_n$, the $q_n$, $q_n + 2$ is a pair of the prime twins.

**Definition 3.** The two-dimensional array $(q_n,\ q_n+2)$ which according to (2.8) ,it be called as **the pseudo prime twins**.

It is now that, the question is whether those pseudo prime twins must be of the prime twins? If the answer is yes, the proof of Lemma 3 is the natural, known as the twin prime conjecture. The next section will explore this issue.

## § 3 The main results and Its prove

Firstly, we give a more general result than theorem 2 in the following.

**Theorem 3.** Considering the sequence $\left\{S_{p_n}\right\}_{n=1}^{\infty}$ in formula（2.3）, For any natural number M , Let $Q(\mathrm{M})$ be the intersection between the collection $\bigcap_{n=2}^{\mathrm{M}} S_{p_n}$ and the interval $[p_{\mathrm{M}},2p_{\mathrm{M}}]$. That is

$$Q(\mathrm{M})=\left(\bigcap_{n=2}^{\mathrm{M}} S_{p_n}\right)\bigcap[p_{\mathrm{M}},2p_{\mathrm{M}}]. \qquad (3.1)$$

If $Q(\mathrm{M})\neq\Phi$ is true, then for $\forall q\in Q(\mathrm{M})$, $q$, $q+2$ will be a pair of the prime twins.

**Proof.** According to (3.1) by the number set $Q(\mathrm{M})$, If $Q(\mathrm{M})\neq\Phi$ , then in case of $q\in Q(\mathrm{M})$, there is $q\in \mathrm{S}_{p_n}$ $(n=2,3,\cdots,\mathrm{M})$ , and $q\in[p_{\mathrm{M}},2p_{\mathrm{M}}]$, thus

$$(q,\ p_n)=1,\ (q+2,\ p_n)=1,\quad p_{\mathrm{M}}<q<2p_M\,(n=2,3,\cdots,\mathrm{M}), \qquad (3.2)$$

By Lemma 2, $q$, $q+2$ will be a pair of the prime twins. This is the proof of Theorem 3.

It is worth pointing out, the set $Q(\mathrm{M})$, by (3.1) describe, is the number set of $q$, the twin prime number $q$, $q+2$ is all meet right on $[p_{\mathrm{M}},2p_{\mathrm{M}}]$. here the precondition is $Q(\mathrm{M})\neq\Phi$ of course. Focus now on whether $Q(\mathrm{M})\neq\Phi$ is set up to each natural number. To answer the questions, we give the simulation experiments first and verify operation by the following Mathematica program:

```
For[n = 2, n <= 2000, n++, q[n_] := {   };
```

```
For[m = (Prime[n] - 1)/2, m <= Prime[n] - 3/2, m++,
If[PrimeQ[2m + 1] && PrimeQ[2m + 3] && 2m + 3 <= 2Prime[n],
q[n] = Union[q[n], {{2m + 1, 2m + 3}}]]]]
Table[q[n], {n, 2, 10}]
Table[q[n], {n, 2, 2000}];
Table[Length[q[n]], {n, 2, 2000}]
```

Computing results are obtained as below:

Out[2]=

{{{3, 5}}, {{5, 7}}, {{11, 13}}, {{11, 13}, {17, 19}}, {{17, 19}}, {{17, 19}, {29, 31}}, {{29, 31}}, {{29, 31}, {41, 43}}, {{29, 31}, {41, 43}}}

Out[4]=（**See Appendix** ）.

The first output is: when the original nine odd prime number *p* (*p* = 3,5,7,11,13,17,19,23,29,31), between prime numbers *p* and 2*p* all twin primes the situation; the second output is: the twin prime number of statistics between prime number *p* and 2*p* for the 1999 smaller odd primes.

Verify operation showed that, as the prime number *p* increased gradually, between prime numbers *p* and 2*p* twin prime number on the overall trend is gradually increasing, although it is not a monotone increasing (see appendix note 3 in the chart). But if between adjacent primes appears to reduce the most only one unit (such as $Q(1999)$ has 235 twin prime number, and there are 234 twin prime number in $Q(2000)$, etc.).Therefore, it is reasonable to speculate that of the prime twins between prime numbers *p* and 2*p* is always present, but also for large prime number *p*, relatively great number of these prime.

Let us discuss the question of whether $Q(\mathrm{M}) \neq \Phi$ is true to each natural number M now.

For every $Q(\mathrm{M})$ defining by (3.1), considering the order of collection, $|Q(\mathrm{M})|$ is a number that $Q(\mathrm{M})$ contains the number of elements. Through the actual calculation can know, when the prime number *p* increases, the prime numbers *p* and $2p$ between the twin prime number for the overall trend is gradually enlarged, appropriate big always seem to M more than 10, such as, $|Q(80)|=12$ , $|Q(2000)|=234$,and so on.

For the convenience of discussion, we put a number $a$ , more than 10 , being called "$a$ is far greater than 1", and to be represented as $a>>1$. So for every natural number $\mathrm{M}$ between 80 ~ 2000, all have $|Q(\mathrm{M})|>>1$.

Let see the below important conclusion.

**Theorem 4.** For large enough the natural number $\mathrm{M}$ , there is always $Q(\mathrm{M})\neq\Phi$ .

**Proof.** For natural numbers $\mathrm{M}$, we make the proof of $|Q(\mathrm{M})|>>1$ by using mathematical induction right now.

First of all, we have seen that $|Q(n)|>>1$ is true for every natural number $n$ between 80 ~ 2000.

Secondly, the inductive hypothesis: For the natural numbers $\mathrm{M}$ （Consider that it is greater than 2000）, there is always

$$|Q(\mathrm{M})|>>1, \quad Q(\mathrm{M})=\left(\bigcap_{n=2}^{\mathrm{M}}S_{p_n}\right)\bigcap[p_{\mathrm{M}},2p_{\mathrm{M}}]. \qquad (3.3)$$

Now lets investigate of natural number $\mathrm{M}+1$，Whether is $|Q(\mathrm{M}+1)|>>1$ or not. Here

$$Q(\mathrm{M}+1)=\left(\bigcap_{n=2}^{\mathrm{M}+1}S_{p_n}\right)\bigcap[p_{\mathrm{M}+1},2p_{\mathrm{M}+1}], \qquad (3.4)$$

On account of

$$\bigcap_{n=2}^{\mathrm{M}+1}S_{p_n}=\left(\bigcap_{n=2}^{\mathrm{M}}S_{p_n}\right)\bigcap \mathrm{S}_{\mathrm{p}_{\mathrm{M}+1}}\subseteq\bigcap_{n=2}^{\mathrm{M}}S_{p_n},$$

$$[p_{\mathrm{M}},2p_{\mathrm{M}}]\bigcap[p_{\mathrm{M}+1},2p_{\mathrm{M}+1}]=[p_{\mathrm{M}+1},2p_{\mathrm{M}}],$$

Using Lemma 2,we have the following conclusions：

（1）For every $q\in S_{p_{\mathrm{M}}}\bigcap[p_{\mathrm{M}+1},2p_{\mathrm{M}}]$ (Note $q\in Q(\mathrm{M})!$), If $q$，$q+2$ are a pair of the twin prime number, then $q\in\bigcap_{n=2}^{\mathrm{M}}S_{p_n}$ . And so $\in S_{p_{\mathrm{M}+1}}$ , that is $q\in\bigcap_{n=2}^{\mathrm{M}+1}S_{p_n}$ ,and then $q\in Q(\mathrm{M}+1)$；

（2）There may exist a pair of the twin prime number $q$， $q+2$ within the interval $[2p_{\mathrm{M}}, 2p_{\mathrm{M}+1}]$. In this case, there will be $q \in \bigcap_{n=2}^{\mathrm{M}+1} S_{p_n}$ , and then $q \in Q(\mathrm{M}+1)$ . Suppose $x(M)$ is the number of the twin prime on the interval $[2p_{\mathrm{M}}, 2p_{\mathrm{M}+1}]$；

（3）Since $p_{\mathrm{M}}, p_{\mathrm{M}+1}$ are adjacent to prime， we can seen that if $q \in S_{p_{\mathrm{M}}} \bigcap [p_{\mathrm{M}}, 2p_{\mathrm{M}}]$ (In this case, $q$， $q+2$ are a pair of the twin prime number!) and $q \notin S_{p_{\mathrm{M}+1}}$ however，it must be $q = p_{\mathrm{M}+1} \in S_{p_{\mathrm{M}}}$ ，but $p_{\mathrm{M}+1} \notin S_{p_{\mathrm{M}+1}}$. In this case, $p_{M+1}$， $p_{M+1}+2$ $q+2$ are a pair of the twin prime number，that is $q = p_{\mathrm{M}+1} \in Q(\mathrm{M})$ and $q = p_{\mathrm{M}+1} \notin Q(\mathrm{M}+1)$.

Based on the above analysis shows that, in terms of the number of $|Q(\mathrm{M}+1)|$ and $|Q(\mathrm{M})|$ relationships are possible results as follows:

$$|Q(\mathrm{M}+1)| = \begin{cases} |Q(\mathrm{M})| - 1, & \text{case1} \\ |Q(\mathrm{M})| + x(M) - 1, & \text{case2} \\ |Q(\mathrm{M})| + x(M), & \text{case3} \end{cases} \qquad (3.5)$$

Here, the three cases in (3.5) corresponding to the three conditions below:

（1） $p_{M+1}$， $p_{M+1}+2$ are a pair of the twin prime number , and $x(M)=0$；

（2） $p_{M+1}$， $p_{M+1}+2$ are a pair of the twin prime number , and $x(M)>0$；

（3） $p_{M+1}$， $p_{M+1}+2$ aren't a pair of the twin prime number , and $x(M) \geq 0$ (Note: In this case, $p_{M+1}+2 \notin \mathrm{P}$, that is $p_{M+1}+2$ being a composite number!).

In either case, the following conclusions are always true.

$$|Q(\mathrm{M}+1)| \geq |Q(\mathrm{M})| - 1. \qquad (3.6)$$

Notice the inductive hypothesis (3.3)，that is $|Q(\mathrm{M})| >> 1$，From （3.3）and （3.6） we have $|Q(\mathrm{M}+1)| \geq |Q(\mathrm{M})| - 1 >> 1$.

In this way, using the mathematical induction principle that, $|Q(\mathrm{M})| >> 1$ is true for these natural numbers M greater than 2000. Therefore $Q(\mathrm{M}) \neq \Phi$ is true in case of M being big enough. This is the proof of Theorem 4.

As a direct corollary of theorem 4, we have the following important conclusion:

**Theorem 5.** For big enough prime $p$ , There always exist a pair of the twin prime number $q$, $q+2$ in the interval $(p, 2p)$.

**Proof.** For being given big enough prime $p$ , assume that it be the $m$-$th$ prime, that is $p = p_m$. In accordance with $\{S_{p_n}\}_{n=2}^{\infty}$ and Theorem 4 as well formula (3.1), we have $Q(m) \neq \Phi$ , and so $q = \min \bigcap_{n=2}^{m} S_{p_n} \in (p, 2p)$. By Theorem 2, $q$, $q+2$ are a pair of the twin prime number in the interval $(p, 2p)$. This is the proof of Theorem 5.

Theorem 5 show that, when the prime number $p$ is large, of the prime twins between prime numbers $p$ and $2p$ is always present. Therefore we have the conclusion that: (2.8) of the definition of pseudo prime twins two-dimensional array is of the prime twins to appropriate large $n$. In fact, we validated the top 2000 primes, in addition to the prime twins does not exist between 5 and 10 foreign (note: 5, 7 is a prime twins, but may appear on the left endpoint of the interval $(5,10)$!), the other between prime number $p$ and $2p$ for there is always the prime twins.

Because there are an unlimited number of prime numbers, so of the prime twins have infinitely many groups. This suggests that the twin prime conjecture that lemma 3 is correct. In fact, we have

**Proof of Lemma 3.** Consider determined by (2.1) and (2.3) sequence $\{S_{p_n}\}_{n=1}^{\infty}$, and assume $n >> 1$.

Based on Theorem 2, $a_n = \min \bigcap_{j=1}^{n} S_{p_j} \in Q(n)$ satisfies: $a_n, a_n + 2$ are a pair of the twin prime number, where $10 < p_n < a_n < 2p_n$. It is obvious that $\lim_{n\to\infty} a_n = \infty$, and so, these twin primes $\{(a_n, a_n + 2)\}_{n=2000}^{\infty}$ have infinitely many groups. This is the proof of Lemma 3.

## §4 Remarks

Through above discussion, we can see the twin prime conjecture is correct, or there are infinitely many group of twin prime number. Careful study theorem 2 ~ 5 proof process, also can get many interesting results and even people who are seeking to prove important conjecture. As a conclusion, we summarized as follows:

（1） Proved the important link of the twin prime conjecture is constructed by the formula (2.1) of the definition of data set, this number is set columns has good properties, the minimum number $q$ that is first $n$ number set intersection can decided a twin prime $q, q+2$.

（2） Proof of theorem 5 process show that each prime number can identify at least one group of twin prime number $q$, $q+2$, this set of twin prime number between $p$ and $2p$. And so easy to prove that the following general results.

**Theorem 6.** When positive $x$ is suitably large, there is at least a pair of twin prime number within interval $(x, 2x)$.

**Proof.** For the appropriate positive number $x$， Let $p(x) = \max\{ q \in P \mid q \le x \}$, then $p(x)$ is the largest prime number less than $x$， and also $p(x) \le x < 2p(x) \le 2x$. By theorem 5，we can assume that $q, q+2$ be a pair of twin prime number within $p(x)$ and $2\,p(x)$. Therefore, $q$, $q+2$ must be a pair of twin prime number within interval $(x, 2x)$. This is the proof of Theorem 6.

（3） By the way, in this paper, the method of tectonic sequence number set (2.3) is actually a new sieve method, the sieve method is never seen in the literature. The author try to use this method to establish a seems to confirm that the famous Goldbach conjecture [3] the sequence number set, of course, how its accuracy to the test of time [25].

Number theory has rich the research content, prime in communication, cryptography and other modern science has a very important application. But there are still many unsolved problems in number theory, including the Fermat number of primes, 1 + 1, etc. There is no doubt that as people the in-depth understanding of the natural numbers, especially the prime rule and found that will inevitably promote the development of number theory science, and the application of number theory in many fields.

**Postscript**: Preliminary work in this paper the source and the author's work in May 2012, was a collection of tectonic sequence $\left\{S_{p_n}\right\}_{n=1}^{\infty}$ to $t_n = \bigcap_{j=1}^{n} S_{p_j}$ in this article, but they were trying to show that $a_n = \min t_n$ satisfy: $a_n, a_n + 2$ prime twins will be failed to settle in $t_n$, $p_n$, $p_n + 2$ prime twins for the estimation of the number of the prime twins. This paper is in my thinking and on June 6, 2013 in 2012[25], based on the structure innovation (2.1) after the formation of the results.

**Appendix**: the first 1999 odd prime, prime numbers $p$ and $2p$ between the twin prime number of statistics:

{1, 1, 1, 2, 1, 2, 1, 2, 2, 2, 3, 3, 2, 2, 3, 4, 3, 3, 4, 3, 4, 4, 4, 6, 7, 6, 6, 5, 5, 7, 7, 8, 7, 8, 7, 8, 8, 8, 8, 9, 8, 8, 7, 7, 6, 7, 8, 8, 7, 8, 8, 7, 7, 7, 8, 8, 7, 7, 7, 6, 7, 8, 9, 8, 8, 10, 10, 10, 9, 9, 9, 9, 9, 9, 9, 9, 9, 9, 10, 12, 11, 12, 11, 11, 12, 12, 12, 12, 11, 11, 11, 11, 11, 11, 11, 11, 13, 12, 14, 15, 15, 15, 15, 14, 15, 15, 15, 15, 14, 14, 14, 15, 14, 14, 15, 14, 15, 16, 16, 16, 16, 16, 16, 16, 16, 16, 17, 18, 18, 18, 19, 20, 20, 20, 20, 20, 20, 20, 21, 21, 21, 20, 20, 19, 20, 21, 21, 20, 21, 21, 21, 20, 20, 21, 21, 21, 21, 22, 23, 23, 23, 24, 24, 25, 25, 25, 25, 26, 26, 27, 26, 26, 25, 25, 27, 26, 27, 26, 27, 28, 28, 27, 27, 27, 27, 27, 28, 28, 29, 28, 29, 30, 30, 30, 31, 31, 31, 31, 31, 31, 30, 30, 30, 30, 31, 30, 30, 30, 29, 30,30, 29, 29, 29, 28, 28, 31, 32, 32, 32, 33, 34, 34, 34, 33, 33, 33, 33, 33, 32, 32, 32, 32, 31, 32, 31, 31, 31, 32, 32, 32, 32, 32, 32, 32, 33, 33, 33, 33, 34, 34, 34, 33, 33, 33, 32, 33, 34, 35, 35, 36, 35, 37, 38, 37, 37, 37, 36, 37, 38, 38, 38, 38, 40, 41, 41, 40, 41, 41, 41, 41, 42, 42, 42, 42, 41, 41, 40, 41, 41, 41, 42, 43, 42, 42, 41, 43, 43, 43, 43, 43, 42, 43, 44, 44, 45, 44, 44, 45, 45, 46, 47, 46, 46, 45, 45, 46, 45, 47, 47, 48, 48, 47, 47, 47, 48, 48, 48, 49, 49, 49, 48, 49, 49, 50, 49, 49, 50, 50, 50, 50, 50, 49, 51, 51, 50, 50, 50, 50, 51, 51, 51, 50, 50, 50, 51, 52, 52, 52, 52, 52, 52, 52, 53, 53, 53, 54, 56, 56, 56, 56, 56, 56, 56, 56, 56, 55, 55, 56, 56, 56, 57, 57, 56, 56, 56, 56, 56, 56, 55, 55, 55, 55, 56, 55, 55, 56, 55, 56, 56, 57, 58, 58, 58, 57, 57, 57, 56, 56, 59, 59, 59, 59, 59, 59, 60, 60, 60, 60, 60, 60, 61, 62, 63, 63, 63, 63, 62, 62, 61, 61, 61, 61, 61, 61, 62, 62, 63,63, 63, 63, 64, 64, 63, 64, 65, 65, 64, 65, 65, 65, 65, 65, 65, 65, 66, 66,65, 65, 64, 64, 66, 65, 65, 65, 65, 65, 65, 65, 65, 66, 67, 66, 66, 65, 66, 66, 67, 67, 68, 69, 69, 69, 68, 68, 67, 69, 69, 69, 69, 69, 68, 68, 68, 67, 67, 67, 66, 67, 67, 66, 66, 67, 67, 67, 67, 67, 67, 67, 68, 69, 68, 69, 69,69, 69, 69, 69, 69,

70, 70, 71, 71, 70, 71, 72, 73, 73, 73, 72, 72, 72, 72, 71, 71, 71, 72, 72, 72, 72, 72, 71, 71, 71, 70,
71, 71, 71, 72, 72, 71, 72,72, 72, 71, 71, 72, 71, 71, 71, 71, 71, 70, 70, 71, 72, 71, 71, 71, 72, 72, \
71, 71, 72, 72, 73, 72, 72, 71, 71, 70, 70, 70, 69, 70, 69, 69, 69, 69, 71, 71, 70, 70, 70, 70, 70, 70,
70, 70, 72, 71, 72, 72, 72, 72, 72, 72, 71, 72, 74, 74, 74, 73, 74, 74, 73, 73, 73, 73, 73, 73, 73, 74,
74, 73, 74, 74, 73, 73, 73, 74, 74, 74, 74, 77, 76, 76, 77, 77, 77, 77, 77, 76, 76, 76, 75, 75,76, 76,
78, 78, 78, 79, 79, 79, 79, 80, 79, 79, 79, 79, 79, 80, 79, 79, 79,79, 79, 79, 80, 79, 80, 79, 80, 81,
81, 82, 82, 82, 82, 81, 81, 81, 81, 82, 83, 84, 84, 84, 84, 84, 84, 85, 86, 85, 85, 86, 87, 87, 86, 86,
86, 86, 86, 86, 86, 86, 87, 87, 87, 87, 87, 87, 87, 86, 87, 87, 87, 86, 87, 88, 88, 87, 87, 87, 86, 86,
86, 85, 85, 86, 87, 88, 88, 88, 89, 90, 90, 90, 89, 89, 89, 88, 88, 87, 87, 88, 88, 88, 88, 88, 88, 88,
88, 87, 88, 89, 89, 89, 89, 89, 89, 89, 89, 89, 89, 89, 89, 89, 90, 90, 89, 89, 88, 89, 89, 90, 90, 90,
90,91, 92, 92, 92, 93, 94, 94, 94, 94, 95, 95, 95, 96, 95, 95, 95, 96, 97, 96,96, 96, 96, 96, 97, 96, 96,
96, 96, 96, 96, 96, 96, 96, 96, 96, 96, 96, 96,95, 96, 96, 96, 96, 96, 96, 96, 96, 95, 95, 95, 95, 95,
95, 96, 96, 96, 95, 96, 96, 96, 96, 98, 98, 98, 98, 97, 97, 97, 96, 96, 96, 96, 96, 97, 97, 97, 97, 96,
97, 97, 97, 96, 97, 96, 96, 96, 96, 96, 96, 95, 95, 94, 94, 93, 93, 93, 93, 92, 92, 93, 95, 96, 96, 95,
96, 96, 96, 96, 97, 98, 97, 98, 97, 98, 98, 98, 98, 98, 98, 99, 100, 100, 100, 100, 101, 101, 101, 101,
101, 101, 101, 102, 101, 101, 101, 102, 102, 102, 103, 103, 102, 102, 103, 103, 103, 103, 103, 105,
106, 106, 105, 106, 106, 105, 105, 104, 104, 104, 104, 104, 104, 105, 105, 104, 104, 104, 104, 103,
103, 103, 103, 103, 103, 103, 103, 103, 102, 102, 101, 102, 102, 102, 102, 101, 101, 101, 101, 102,
102, 103, 104, 104, 105, 105, 105, 105, 105, 105, 105, 105, 105, 105, 105, 104, 104,105, 105, 106,
107, 107, 107, 108, 109, 109, 108, 108, 108, 108, 108, 108, 108, 108, 109, 108, 108, 109, 109, 108,
108, 110, 110, 110, 110, 111, 111, 110, 110, 111, 111, 112, 112, 112, 112, 112, 112, 112, 113, 113,
113, 112, 113, 112, 112, 112, 112, 112, 112, 112, 112, 111, 111, 111, 112, 113, 114, 114, 114, 114,
114, 113, 114, 114, 114, 113, 113, 113, 114, 114, 115, 115, 116, 116, 116, 115, 115, 115, 115, 115,
116, 115, 116, 116, 116, 115, 115, 116, 116, 116, 116, 116, 116, 116, 117, 117, 118, 118, 118, 118,
118, 119,119, 119, 119, 119, 121, 121, 121, 120, 121, 121, 120, 121, 121, 120, 120, 121, 121, 123,
123, 123, 123, 123, 125, 125, 124, 126, 125, 125, 125, 124,126, 127, 126, 126, 126, 128, 128, 128,
128, 129, 129, 129, 130, 131, 131,131, 131, 131, 131, 131, 131, 131, 131, 131, 130, 130, 132, 132,
131, 131, 131, 131, 131, 131, 131, 130, 130, 130, 130, 130, 130, 130, 130, 130, 129, 129, 128, 128,
127, 129, 128, 128, 128, 128, 128, 128, 128, 128, 128, 128, 129, 129, 130, 131, 132, 132, 132, 132,
131, 131, 131, 131, 131, 130, 130,131, 133, 132, 132, 133, 133, 133, 133, 132, 133, 133, 133, 133,
133, 133, 133, 133, 133, 134, 134, 133, 133, 134, 134, 134, 134, 135, 135, 134, 134, 135, 135, 135,
137, 136, 137, 136, 136, 136, 135, 136, 136, 135, 135, 135, 135, 136, 136, 135, 135, 135, 135, 135,
135, 136, 136, 136, 137, 138, 138, 138, 139, 139, 139, 138, 139, 139, 138, 138, 139, 139, 138, 138,
138, 138, 139, 141, 141, 142, 141, 141, 142, 142, 141, 141, 141, 141, 142, 141, 143, 143, 143, 143,
143, 143, 144, 144, 144, 144, 144, 144, 144, 144, 144, 145,145, 145, 146, 146, 145, 145, 145, 145,
145, 146, 147, 148, 148, 149, 152, 152, 152, 152, 152, 151, 151, 152, 152, 151, 151, 151, 152, 151,
151, 151, 151, 151, 151, 151, 151, 152, 153, 154, 153, 153, 152, 153, 153, 153, 153, 153, 152, 152,
154, 154, 153, 153, 152, 152, 153, 153, 153, 154, 154, 154, 154, 155, 155, 156, 156, 157, 157, 158,
158, 159, 158, 158, 159, 159, 159,159, 159, 160, 160, 160, 160, 160, 161, 161, 160, 160, 160, 161,
162, 163, 162, 162, 162, 162, 162, 163, 163, 163, 164, 164, 164, 165, 165, 164, 164, 163, 163, 163,
164, 163, 164, 164, 164, 164, 164, 165, 165, 165, 164, 165, 166, 166, 167, 167, 167, 167, 168, 168,
168, 168, 167, 167, 168, 168, 167, 167, 167, 167, 167, 167, 167, 166, 166, 167, 166, 167, 167, 167,
166, 166, 166, 166, 166, 166, 166, 165, 166, 166, 167, 167, 167, 166, 166, 165, 165, 165, 165, 165,
165, 165, 165, 165, 165, 165, 165, 165, 164, 164, 164, 164, 164, 164, 164, 164, 164, 164, 165, 165,

165, 166, 166, 166, 166, 167, 167, 167, 166, 166, 166, 166, 166, 166, 167, 167, 167, 166, 166, 166, 166, 166, 167, 168, 168, 168, 168, 168, 169, 169, 170, 170, 170, 170, 170, 171, 171,172, 172, 171, 171, 171, 171, 171, 171, 171, 172, 172, 172, 171, 171, 172, 172, 172, 173, 173, 173, 173, 174, 173, 173, 172, 172, 172, 172, 172, 173, 173, 173, 173, 173, 173, 174, 175, 175, 175, 175, 175, 175, 175, 175, 175,174, 174, 174, 174, 174, 174, 174, 174, 174, 174, 174, 174, 174, 173, 177, 177, 177, 176, 176, 176, 176, 178, 179, 179, 179, 179, 180, 180, 180, 180, 180, 181, 181, 182, 182, 182, 182, 182, 182, 183, 183, 184, 184, 184, 183,183, 183, 182, 182, 183, 182, 182, 181, 181, 182, 182, 181, 181, 183, 183,184, 184, 184, 183, 183, 184, 185, 186, 185, 185, 186, 185, 185, 185, 185, 185, 184, 185, 185, 186, 185, 185, 184, 184, 184, 185, 186, 186, 186, 185,185, 186, 187, 187, 187, 188, 188, 189, 189, 190, 190, 190, 190, 189, 190, 191, 191, 192, 191, 191, 192, 192, 192, 193, 192, 192, 192, 192, 192, 192, 192, 192, 192, 191, 191, 191, 191, 191, 192, 192, 192, 192, 192, 191, 191, 191, 190, 191, 190, 191, 191, 190, 190, 190, 190, 190, 190, 190, 191, 192, 192, 192, 192, 192, 192, 192, 192, 192, 192, 192, 192, 192, 193, 193, 193, 193, 193, 194, 194, 194, 193, 193, 194, 194, 194, 194, 194, 194, 195, 195, 195, 195, 195, 196, 196, 196, 197, 197, 198, 198, 198, 198, 198, 198, 198, 198, 199, 198, 198, 198, 198, 198, 198, 199, 199, 199, 199, 200, 201, 201, 201, 200, 200, 201, 200, 200, 200, 200, 200, 200, 199, 199, 199, 198, 198, 198, 198, 198, 198, 198, 200, 201, 201, 201, 201, 201, 201, 201, 201, 201, 201, 202, 202, 202, 203, 204, 203, 204, 204, 204, 205, 205, 204, 204, 203,204, 204, 204, 204, 204, 205, 205, 204, 204, 203, 203, 204, 204, 205, 205,205, 205, 205, 205, 205, 205, 205, 207, 207, 208, 207, 207, 207, 207, 207, 207, 208, 208, 208, 207, 207, 207, 207, 209, 209, 210, 209, 209, 208, 209, 209, 209, 210, 210, 210, 210, 210, 209, 211, 212, 211, 211, 212, 213, 213, 212, 212, 212, 213, 213, 214, 215, 215, 215, 215, 216, 215, 215, 215, 216, 217, 217, 217, 217, 217, 217, 216, 218, 218, 219, 219, 219, 219, 220, 220, 220, 220, 220, 222, 222, 222, 222, 221, 222, 221, 221, 221, 222, 223, 222, 222, 222, 222, 222, 222, 222, 222, 222, 225, 225, 225, 224, 224, 224, 225, 225, 226, 226, 225, 227, 227, 227, 227, 227, 227, 227, 226, 226, 226, 226, 227, 227, 226, 226, 226, 226, 226, 227, 228, 228, 229, 229, 229, 230, 231, 231, 231, 232, 231, 231, 231, 230, 230, 231, 233, 233, 232, 233, 233, 233, 234, 234, 234, 234, 234, 234, 235, 235, 234}

Note **1**: The number $n$ in table primes between $p_{n+1}$ and $2p_{n+1}$ twin logarithm of prime numbers.

Note **2**: One exception in the above table, that is the second number "1" in the table just corresponding primes of (5, 7), it is the twin prime number on $[5, 10]$, but it is not the twin prime number within （5，10）.

Note **3**: The number of $n$ in the above table with the corresponding point series $\{(n, p_{n+1})\}_{n=1}^{1999}$ of prime $p_{n+1}$ relationship here is as follows:

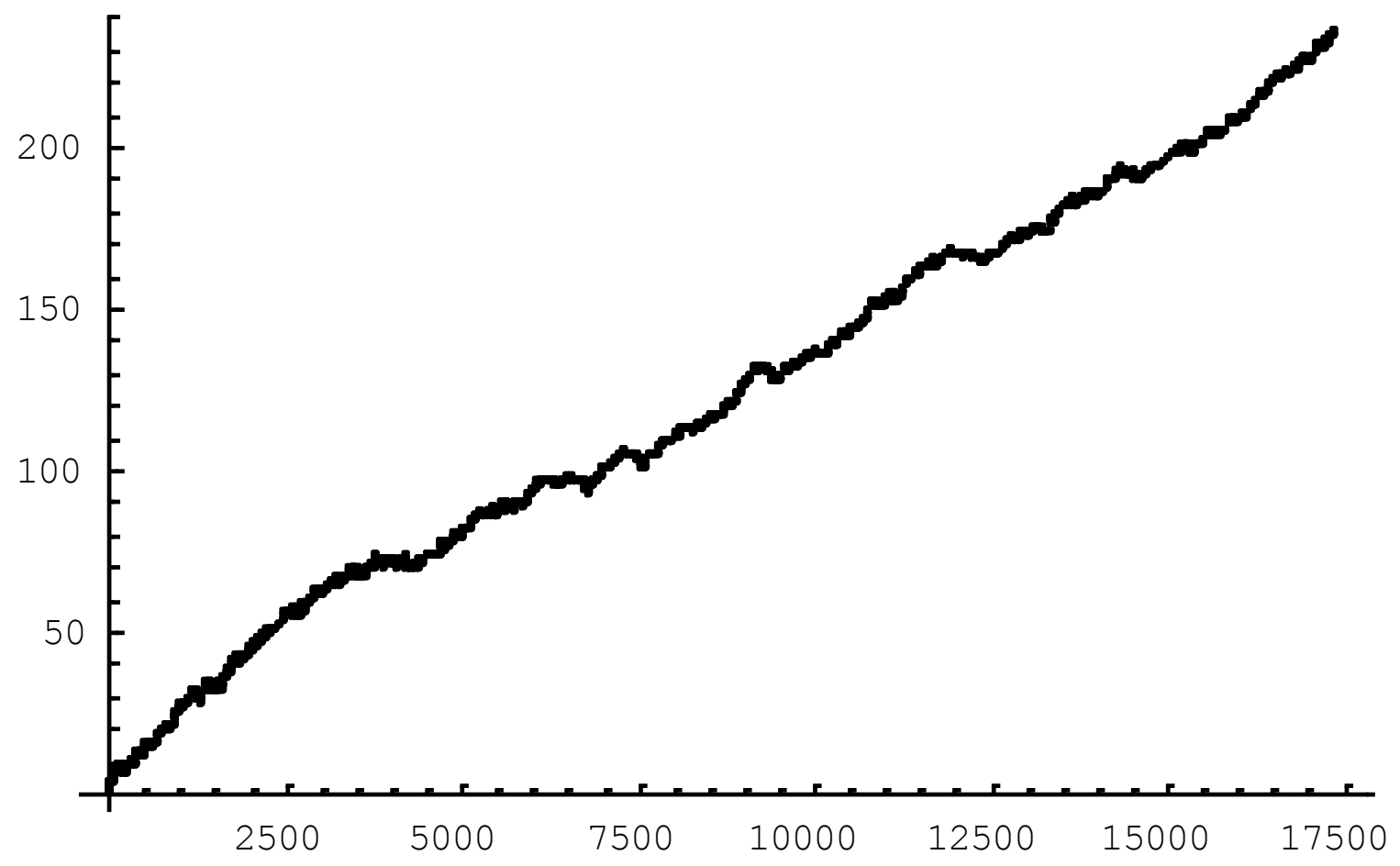

**References**

[1] Hua Luogeng waiting. Number theory guidance [M], science press，1959，1~45.

[2] Yu.I.Manin & A.A.Panchishkin. Introduction to Modern Number Theory(second Edition)[M]. Springer Verlag Berlin Heidelberg. 2005.

[3] Richard E. Crandall, Carl Pomerance. "Prime numbers: a computational perspective". Springer. 2005: pp.14.

[4] Hector Zenil. Goldbach's Conjecture , Wolfram Demonstrations Project, 2007.

[5] Eric W. Weisstein, Goldbach Number, Math World.

[6] Wang Yuan. The Goldbach Conjecture. World Scientific Publishing Company (second Edition)[M]. 2002.1~18.

[7] Hardy, G. H. and Littlewood, J. E.. Some Problems of Partitio Numerorum (III): On the expression of a number as a sum of primes. Acta Mathematica. 1923(44): 1~70.

[8] Pan Chengdong and Pan Chengbiao. Goldbach conjecture (first edition) [M]. Science press. 1981. 148~149.

[9] J M Deshouillers, G Effinger, H Te Riele, D Zinoviev. A Complete Vinogradov 3-Prime Theorem under the Riemann Hypothesis. Electronic Research Announcements of The American Mathematical Society, 3: 99-104

[10] J.-M. Deshouillers, H. J. J. te Riele and Y. Saouter. New experimental results concerning the Goldbach conjecture. Lecture Notes in Computer Science.1423/1998: 204~215

[11] J. R. Chen, On the representation of a larger even integer as the sum of a prime and the product of at most two primes. Sci. Sinica 16 (1973), 157~176.

[12] Pipping,Nils. "Die Goldbachsche Vermutung und der Goldbach-Vinogradovsche Satz." Acta. Acad. Aboensis, Math. Phys. 1938(11): 4~25.

[13] A. Granville, J. van de Lune, and H. J. J. te Riele. Checking the Goldbach conjecture on a vector computer. Number Theory and Applications, R. A. Mollin (ed.) Kluwer

Academic Press. 1989: 423~433.

[14] Matti K. Sinisalo, Checking the Goldbach conjecture up to $4 \cdot 10^{11}$, Mathematics of Computation, vol. 61, 1993 (204): 931~934.

[15] Jörg Richstein, Verifying the Goldbach conjecture up to $4 \cdot 10^{14}$, Mathematics of Computation, vol. 70, 2000(236):1745~1749.

[16] Tomás Oliveira e Silva. Goldbach conjecture verification , 2012.

[17] G. J. Rieger. Solution of the Waring-Goldbach problem for algebraic number fields. Bull. Amer. Math. Soc., 68, 1962 (3): 234~236.

[18] Richard E. Crandall, Carl Pomerance. "Prime numbers: a computational perspective". Springer. 2005: pp.14.

[19] Xinhuanet. How long can Goldbach conjecture " conjecture ", xinhua news agency, August 20, 2002.

[20] Tian. Song. Theory of civil science fans why can't a scientific sense of success?[J]. Science, technology and dialectics. 2004(3): 108~112.

[21] The China youth daily. Chinese decipher "twin prime conjecture" influence or the super trained Chen jingrun 1 + 2, http://news.ifeng.com/world/detail_2013_05/18/25438866_0.shtml?_newshao123，18 May 2013

[22] Beijing university news. Mathematical sciences alumni Zhang Yitang major breakthroughs were made in prime twins in the research, http://pkunews.pku.edu.cn/xwzh/2013-05/15/content_272348.htm，15 May 2013

[23] Nature. First proof that infinitely many prime numbers come in pairs，http://www.nature.com/news/first-proof-that-infinitely-many-prime-numbers-come-in-pairs-1.12989，14 May 2013

[24] Davide Castelvecchi. Mathematicians come closer to solving Goldbach's weak conjecture. http://www.nature.com/news/mathematicians-come-closer-to-solving-goldbach-s-weak-conjecture-1.10636

[25] Baoshan Zhang. Matrix Online Characteristic Number and Goldbach’s Conjecture. http://arxiv.org/abs/1306.0795. 4 Jun 2013